\documentclass[12pt,a4paper]{article}
\usepackage[a4paper,margin=2cm]{geometry}

\usepackage{amsmath}
\usepackage{amsthm}
\usepackage{subcaption}
\usepackage{amssymb}

\usepackage{graphicx}
\usepackage{nccmath}
\usepackage{setspace}
\usepackage[colorlinks=true, allcolors=blue]{hyperref}
\usepackage{enumerate}
\usepackage{multirow}
\usepackage{float}
\usepackage[dvipsnames]{xcolor}
\usepackage{setspace}
\usepackage{tikz}
\usetikzlibrary{positioning}
\usepackage{booktabs}
\usepackage{enumitem}

\newtheorem{theorem}{Theorem}

\newtheorem{definition}{Definition}

\theoremstyle{definition}
\newtheorem{example}[theorem]{Example}

\def \mod#1{{\:({\rm mod}\ #1)}}

\let\oldproofname=\proofname
\renewcommand{\proofname}{\textup{\textbf{\oldproofname}}}

\title{Poetry of Repetition: Constructing Verse through Combinatorial Design Theory}

\author{Ajani De Vas Gunasekara \footnote{School of Arts and Sciences, The University of Notre Dame Australia, Sydney NSW 2007, Australia} \qquad Miriam Wei Wei Lo \footnote{
School of Humanities and Social Sciences, Sheridan Institute of Higher Education, Perth WA 6000, Australia}
}
\date{}
\begin{document}

\maketitle

 \begin{abstract}
 
  This paper investigates the connections between combinatorial design theory and the creation of new forms of poetry through a specific combinatorial structure called \emph{Steiner triple systems}. We introduce five original poems constructed using variations of Steiner triple systems on seven and nine words, illustrating how mathematical structures can inform and inspire new poetic forms. The work includes a reflective discussion from dual creative perspectives; one emphasizing structural design and the other, literary expression, highlighting how formal constraints can foster occasional frustrations, but also unexpected artistic freedoms. This study demonstrates the potential of mathematical frameworks as generative tools in literary creativity. 

\end{abstract}

\section{Introduction}

\emph{Combinatorics} is a branch of mathematics concerned with counting, arrangement, and selection of objects. \emph{Combinatorial design theory} focuses on the systematic arrangement and grouping of distinct, separate objects according to specific rules and patterns. It explores how these objects can be combined or organized to achieve certain desired properties, often ensuring balance, fairness, or symmetry within the structure. This focus aligns with Mirsky's broader view of combinatorics, who describes it as:\cite{Mirsky1979},

\begin{quote}
  Combinatorics is a range
of linked studies which have something in common and yet diverge widely in
their objectives, their methods, and the degree of coherence they have
attained. Most are concerned with criteria for the existence of certain
`patterns' or `arrangements' or `configurations', where these terms need to be
interpreted in a very broad sense. 
\end{quote}

In this sense, combinatorial design theory can be seen as one particular area within combinatorics that studies specific types of patterns and arrangements, driven by the goal of achieving balance and structure.

Poetry is a literary art form that deploys the aesthetic, figurative, and aural qualities of language to evoke meaning in ways that both incorporate and go beyond literal interpretation. Jill. P. Baumgaertner observes, ``poetry, like all other artistic expression, is an attempt to name that experience, to create feeling, to express the otherwise inexpressible" (p.~2)~\cite{Baumgaertner1990}. Individual works of poetry, or poems, typically incorporate a number of structural constraints. For example, most poems are structured into lines and stanzas (with the notable exception of the prose-poem, see Hetherington and Atherton \cite{HetheringtonAtherton2020}). Many poems also incorporate aural patterns (such as rhyme and metric schemes) as additional structural constraints. These structural constraints have the capacity to intensify expressions of emotion and thought.  As American poet Robert Hass says in relation to the 4/4/4/2 structure of the Shakespearean sonnet: ``say it, say it, contradict or qualify it, nail it" (p.~185)~\cite{Hass2017}. Poetry's relationship with patterns, repetitions, constraints, and combinations lends itself to structural exploration. This is where the field of combinatorics intersects with poetry.

There is a long and rich history connecting poetry and mathematics. Over the centuries, poets and mathematicians alike have been drawn to the interplay between structure and expression, constraints and creativity. A significant number of individuals have explored this intersection, highlighting how mathematical ideas can inform poetic form, and how poetry can reflect mathematical thinking. Notable contributions include the work of Sarah Glaz \cite{Glaz2011}, who has written extensively on the use of mathematical themes and structures in poetry; Daniel May \cite{May2021}, who examines how mathematical models can inspire poetic composition; and Sarah Hart \cite{Hart2023}, whose recent work investigates the broader cultural and aesthetic connections between mathematics and the literary arts.

As Karaali and Lesser note, there are three broad (and often overlapping) categories of mathematical poetry \cite{KaraaliLesser2020}:

\begin{itemize}
    
    \item poems that take mathematics as their primary subject matter,
    \item poems that use mathematical language or imagery to explore non-mathematical topics, and
    \item poems whose structure is shaped or influenced by mathematical principles.
\end{itemize}

The focus of this paper is on the third category (poems whose structure is shaped or influenced by mathematical principles). Combinatorics offers powerful tools both to analyze and construct poetic forms. On the analytical side, it helps uncover the structural patterns that underlie existing poetry such as the intricate word rotation in a \emph{sestina}, where the six end-words are permuted in a fixed sequence across stanzas. On the constructive side, combinatorial techniques can be used to generate new poetic constraints, such as limits on the number of stanzas, the number of words, and the repetition of certain words or lines. While poetry has long been a playground for creativity and logic, embracing and sabotaging patterns, combinatorics highlights and extends this capacity, offering new ways to explore the balance between form and expression.

As early as the first century BCE, poets were already exploring the combinatorial nature of poetic composition. In Book II of ``De rerum natura" (``On the Nature of Things"), the Roman poet and philosopher Lucretius presents poetry itself as an inherently combinatorial act \cite{May2021}.

Before the twentieth century, poetic works often employed combinatorial techniques that offered writers multiple interpretive possibilities; sometimes generating thousands of variations from a single short text. One such example is the \emph{cento}, a poetic form composed entirely of lines borrowed from other poems, rearranged to create new meaning. The structure of a cento reflects combinatorial creativity, as its effect depends on the selection and permutation of pre-existing lines.

In the mid-twentieth century, this interplay became more intentional. The Oulipo group, founded in the 1960s, brought together writers dedicated to experimenting with constraints and mathematical structures in literary creation. It is worth noting that the Oulipo group also wanted to resist expressivism, which means that the experience of reading Oulipo work ``is like the pleasure of watching someone build something rather than the pleasure of having someone tell you something" (Hass, 2017, p.~379)~\cite{Hass2017}.

In \cite{Barriere2017}, L. Barrière has explored the intersection of combinatorics and the arts by posing the foundational question: how can combinatorial structures be used in artistic creation? It surveys twentieth-century works in music, literature, and visual arts that incorporate combinatorial thinking, examining both the presence of such structures and the artists' intentions. The study highlights patterns, differences, and commonalities across these disciplines, offering an initial mapping of combinatorics in creative practice.

In this paper, we explore the connections between combinatorial design theory and the creation of new forms of poetry through a specific combinatorial structure called ``Steiner triple systems".  We present five original poems constructed according to variations of Steiner triple systems on seven and nine words. These poems serve as both artistic expressions and structural demonstrations of the underlying combinatorial principles. Alongside the poems, we provide a reflective discussion on the creative process from two complementary perspectives: one grounded in structural and mathematical reasoning, while the other guided by linguistic intuition and poetic sensibility. We examine how the imposition of combinatorial constraints shaped the selection, arrangement, and repetition of words, and how these limitations paradoxically opened new creative possibilities. Through this dialogue between structure and expression, we aim to highlight the potential of mathematical forms as generative frameworks for literary creativity.

\section{Applying Steiner Triple Systems in poetry  construction}

A \emph{Steiner triple system} of order \(u\), denoted \(\mathrm{STS}(u)\), is a pair \((U,\mathcal{A})\) where \(U\) is a set of \(u\) points and \(\mathcal{A}\) is a collection of three-element subsets of \(U\) (called \emph{blocks} or \emph{triples}) such that every unordered pair of distinct points of \(U\) lies in exactly one triple. Equivalently, one asks for a way to form triples from \(U\) elements so that each pair appears together precisely once.

Steiner triple systems were first defined by  Woolhouse in 1844 in the 1733 prize question of \emph{Lady's and Gentleman's Diary} \cite{ColRos, kirkman1847}. 

\begin{quote}
    ``Determine the number of combinations that can be made out of $n$ symbols, $p$ symbols in each; with this limitation, that no combination of $q$ symbols, which may appear in any one of them  shall be repeated in any other."
\end{quote}

In 1847, Kirkman gave an important answer to the $p = 3$ and $q = 2$ case of this problem by proving that the necessary and sufficient conditions for the existence of a Steiner triple system of order $u$ is that $u$ leaves a remainder of 1 or 3 upon division by 6. Denoted as, $u \equiv 1,3 \mod{6}$. Thus, examples of such values include $ u = 3,7, 9,13, 15, \ldots$.

\begin{definition}[Resolvable Steiner triple systems] 
A resolvable Steiner triple system (STS) is a Steiner triple system in which the collection of triples can be divided into groups called parallel classes. Each parallel class partitions the entire point set, meaning that every point appears in exactly one triple within each class. In other words, the triples can be organized into subsets where each subset covers all points exactly once, and together, these subsets include all triples of the system.
    
\end{definition}

These resolvable Steiner triple systems are also known as \emph{Kirkman triple systems}. Lu \cite{Lu1965} and Ray-Chaudhuri and Wilson \cite{RayChaudhuriWilson1971} independently proved that a resolvable Steiner triple system of order $u$ exists when $u \equiv 3 \pmod{6}$. In other words, $u$ must be a positive integer that leaves a remainder of 3 when divided by 6. Examples of such values include $u = 3, 9, 15, 21, \ldots$.

Another remarkable problem in the history of Steiner triple systems is Kirkman's schoolgirl  problem (1850) \cite{ColRos}.

\begin{quote}
    ``Fifteen young ladies in a school walk out three abreast for seven days in succession: it is required to arrange them daily, so that no two shall walk twice abreast."
\end{quote}

The answer to this problem is a resolvable Steiner triple system of order 15 . This recreational problem got the attention of many, and a number of mathematicians studied the problem and its generalizations.

Steiner triple systems \((\mathrm{STS})\) provide a precise combinatorial framework for poetry. In this study, each point corresponds to a unique word, which we refer to as a \emph{keyword}, and each triple corresponds to a ``line'' of three keywords (possibly with some filler words). The system is arranged so that every unordered pair of distinct keywords appears together in exactly one line. A Steiner triple system of order \(u\) \((\mathrm{STS}(u))\), therefore, contains exactly \(u\) keywords. By definition, a Steiner triple system is a pair \((U, \mathcal{A})\), where \(U\) is the set of \(u\) keywords and \(\mathcal{A}\) is the collection of lines of the poem, such that every pair of keywords appears together exactly once in some line. We refer to a poem following this structure as a Steiner triple poem.

In a Steiner triple poem with \(u\) keywords, each keyword appears in exactly \(\tfrac{u-1}{2}\) lines. This quantity is always an integer because \(u\) must be odd, with \(u \equiv 1 \text{ or } 3 \pmod{6}\). Furthermore, the total number of lines is given by \(\tfrac{u(u-1)}{6}\). 

This construction produces a richly interwoven poem that exhibits balance and symmetry: every keyword recurs exactly \(\tfrac{u-1}{2}\) times, and each pair of keywords meets precisely once, offering a fresh perspective with every repetition. These keywords naturally shape the theme or central subject of the poem.

The use of recurring words repeated according to a formula is not new to poetry. The pantun (pantoum), villanelle, sestina, and ghazal are existing poetic forms that have emerged from diverse cultures; all play with some combination of word and/or line repetition \cite{PoetryFoundationForm}.  What the Steiner triple system offers is a different set of rules for the poetic game of repetition. It is an extension of the game that presents a novel set of constraints and freedoms. 

Imposing these combinatorial constraints inevitably narrows the poet's choice of words, but this very narrowing can be liberating. By working within a fixed set of possibilities, the poet discovers a kind of freedom in limitation, a balance that lies at the heart of poetic composition.

 We explore three variations of Steiner Triple Poems:

\begin{enumerate}
    \item \textit{Pure Steiner Triple Poems} - These contain only the designated keywords, without any additional (filler) words. In this version, the underlying mathematical structure is explicitly visible, and the combinatorial design is clearly traceable through the poem’s construction.

    \item \textit{Relaxed Steiner Triple Poems} - These include filler words or natural language expressions interwoven with the  keywords. While the mathematical design is still preserved, it is less visible, allowing for greater linguistic fluidity and expressive potential.

    \item \textit{Resolvable Steiner Triple Poems} - These poems follow the structure of a resolvable Steiner triple system. Each parallel class of triples functions as a stanza, with every line containing three distinct keywords and each keyword appearing exactly once per stanza. The poems can be either pure - using only the designated keywords - or relaxed, where additional filler words are included for narrative or expressive flexibility.
\end{enumerate}

Each variation serves a distinct purpose: the pure version emphasizes mathematical elegance and structural clarity; the relaxed version prioritizes poetic readability and interpretive richness; and the resolvable version introduces an added layer of organization by grouping lines into stanzas that mirror parallel classes in a resolvable design. The choice among these forms depends on the desired balance between formal combinatorial structure and literary expression.

\subsection{Graph theoretic representation of a Steiner triple poem}

In this section, we show how a Steiner triple poem can be represented using a \emph{graph}. This representation also shows how the lines of the poem are formed, each containing exactly three keywords.

A graph $G = (V,E)$ is a mathematical model for a network, consisting of:
\begin{itemize}
    \item A set of \emph{vertices} (or \emph{nodes}) $V$, which one may visualize as points (balls).
    \item A set of \emph{edges} $E$, each of which joins exactly two vertices and can be visualized as a line between those points (balls).
\end{itemize}
To indicate a relationship or connection between two entities, we can place an edge between their corresponding vertices. This abstraction allows us to represent various real‐world systems in a uniform way.

\begin{definition}[Complete graph]
A \emph{complete graph} is a graph where every pair of distinct vertices is connected by exactly one edge. That means each vertex is directly linked to every other vertex. Note that a complete graph on $n$ vertices is denoted by $K_n$, where $n$ is a positive integer.
\end{definition}

\begin{definition}[Subgraph]
A \emph{subgraph} is a graph formed by selecting a subset of the vertices and edges of a larger graph. The selected edges must connect only the selected vertices and must also exist in the original graph.
\end{definition}

A Steiner triple poem on $u$ keywords can be represented naturally using a graph. In this representation, each keyword corresponds to a vertex, and two vertices are connected by an edge if the corresponding keywords appear together in the same line of the poem. Because every pair of keywords appears together exactly once in a Steiner 
triple system, the resulting graph contains all possible edges between vertices. Hence, the structure of the poem corresponds to a complete graph on \(u\) vertices, denoted \(K_u\).

Starting from the graph representation of a Steiner triple poem, we can recover its lines by applying a graph-theoretic technique called \emph{graph decomposition}. In this setting, a decomposition partitions the edges of the complete graph into disjoint subgraphs, each corresponding to a triple of keywords. This provides a 
systematic way to reconstruct the triples of the Steiner triple system, and hence the lines of the poem.  

\begin{definition}[Graph decomposition]  
A graph decomposition is a partition of the edge set of a graph into disjoint subsets, where each subset induces a smaller subgraph of a specified type.  
\end{definition}  

A central example is a \emph{triangle decomposition}, in which the edges of a graph are grouped into triples, each forming a triangle that appears in the original graph.  
In particular, a triangle decomposition of the complete graph on \(u\) vertices, \(K_u\), covers all edges of \(K_u\) with edge-disjoint triangles, so that every edge belongs to exactly one triangle. Thus, a triangle decomposition of \(K_u\) is equivalent to a Steiner triple system of order \(u\) \((\mathrm{STS}(u))\).  

This equivalence provides a natural mechanism for generating Steiner triple poems. The chosen keywords act as the vertices of a complete graph, and a triangle decomposition specifies how the words are grouped into poetic lines. Each triangle corresponds to a line of the poem, ensuring that every pair of keywords appears together exactly once across the entire composition. In this way, the combinatorial 
structure determines the poem’s framework while still allowing flexibility in its expressive form.

\section{Examples of Steiner Triple Poems}

In this section, we present some original examples of Steiner triple poems of each type for seven and nine keywords. Instead of providing a separate graph representation for each poem, we use two graphs~$K_7$ and $K_9$~to illustrate triangle decompositions corresponding to the respective poems. This is sufficient, as the structure is similar across all variations: only the keywords change, serving as the vertices of the graph.

\subsection{\texorpdfstring{$\mathrm{STS}(7)$}{STS(7)}  - Poems with seven Keywords}

\(\mathrm{STS}(7)\) is a Steiner Triple System of order 7. A poem that conforms to this structure contains seven distinct keywords. Since each pair of keywords must appear together exactly once, the poem consists of $\frac{7(7 - 1)}{6} = 7$ lines, with exactly three keywords per line. Consequently, each keyword appears in precisely three different lines.

An example of a pure Steiner Triple poem with seven keywords, one that uses only the designated keywords without any additional words, is discussed by Sarah Hart in \cite{Hart2023}.

Below, we present our constructions of both a pure Steiner Triple poem and a relaxed version. In the relaxed version, filler words are interwoven with the keywords to allow for greater narrative or expressive flexibility. We also illustrate the underlying graph representations of these poems.

\begin{example}
   Pure Steiner triple poem: 
\\   

\textbf{Karak}*
\\
\\
\indent Bird seeks home.

Seeks trees here.

Bird food here?

Trees: Food, home.

\vspace{0.4 cm}

Trees, Bird? Not …

Seeks food? Not …

Home not here.
\\

*Karak is the Noongar word for a red-tailed black cockatoo.

\vspace{0.4 cm}

\textbf{Discussion of the poem:} This poem, composed by the second author (Lo), began with expressive intent: Lo first considered what she wanted to write about. Lo had recently moved to Hilton, a suburb in Walyalup / Fremantle, Western Australia. Many residents in this suburb are concerned about loss of habitat for the karak (red-tailed black cockatoo) and express their concern by putting signage in their gardens and by campaigning for other residents to preserve or plant trees for karak in their gardens. Lo envisaged this poem as her contribution to the campaign. She chose keywords in keeping with this expressive intent. She then used the triangle decomposition of $K_7$ (complete graph on seven vertices) to help her visualise the word groupings, which she noted down in line form. As line order did not matter, she could rearrange the word groupings to suit intent. She could also organise word order within the line to suit intent. Working from within these constraints, with a strong expressive intent, the poem emerged. Punctuation became especially useful in creating nuances of expression - the use of ellipses (...) in the second stanza, for example, allows the reader to fill in the blanks. Perspective also allowed for some flexibility, with the poem moving from the voice of the karak in the first stanza to the voice of a human in the second stanza.

To visualize the structure graph-theoretically, we represent each keyword as a vertex and construct the complete graph on seven vertices, denoted $K_7$. Each line of the poem corresponds to a triangle in this graph. Therefore, the poem forms a triangle decomposition of $K_7$, in which the seven triangles collectively cover all the edges. Each triangle, representing a line of the poem, contains exactly three keywords.

Keywords (vertex set): \{Bird, Seeks, Home, Trees, Here, Food, Not\}

\begin{figure}[h!]
\centering
\begin{tikzpicture}[scale=3, every node/.style={circle, draw, minimum size=1.75cm}]
  % Define the seven nodes in a circular layout
  \node (Bird)    at (90:1)   {Bird};
  \node (Seeks)     at (40:1)   {Seeks};
  \node (Here)     at (-10:1)  {Here};
  \node (Not)    at (-60:1)  {Not};
  \node (Trees)  at (-110:1) {Trees};
  \node (Home)  at (-160:1) {Home};
  \node (Food)    at (150:1)  {Food};

  % Triangle decomposition: 7 disjoint triangles covering all edges
  \draw[thick, blue]     (Bird) -- (Seeks) -- (Home) -- (Bird) --cycle;
  \draw[thick, red]      (Bird) -- (Food) -- (Here) -- (Bird) -- cycle;
  \draw[thick, green!70!black] (Bird) -- (Not) -- (Trees) -- (Bird) -- cycle;
  \draw[thick, orange]   (Seeks) -- (Not) -- (Food) -- (Seeks) -- cycle;
  \draw[thick, black]   (Seeks) -- (Trees) -- (Here) -- (Seeks) -- cycle;
  \draw[thick, teal]     (Home) -- (Trees) -- (Food) -- (Home) -- cycle;
  \draw[thick, magenta]  (Here) -- (Not) -- (Home) -- (Here) --  cycle;

\end{tikzpicture}
\caption{Triangle decomposition of the complete graph $K_7$ on seven keywords from the poem ``Karak". Each triangle corresponds to one line in the poem.}
\label{fig:sts7-poem1}
\end{figure}
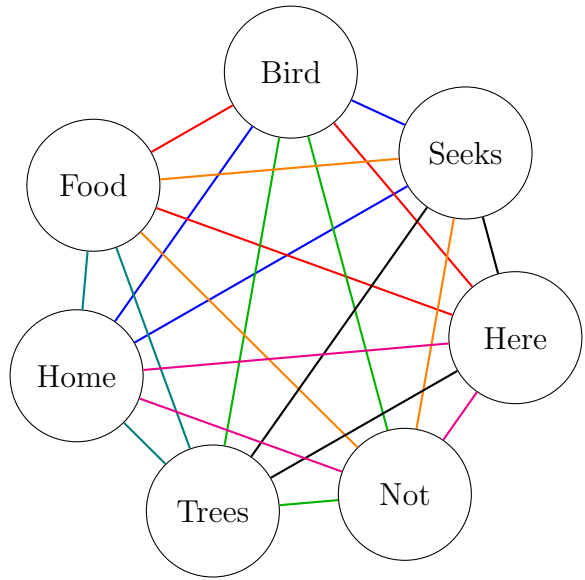

\bigskip

\begin{table}[H]
\centering
\begin{tabular}{@{}lll@{}}
\toprule
\textbf{Line \#} & \textbf{Triangle (Color)} & \textbf{Vertices} \\ \midrule
1 & Blue     & \{Bird, Seeks, Home\} \\
2 & Black    & \{Seeks, Trees, Here\} \\
3 & Red      & \{Bird, Food, Here\} \\
4 & Teal     & \{Home, Food, Trees\} \\
5 & Green    & \{Bird, Not, Trees\} \\
6 & Orange   & \{Seeks, Not, Food\} \\
7 & Magenta  & \{Here, Not, Home\} \\
\bottomrule
\end{tabular}
\caption{The 21 edges of \( K_7 \) are partitioned into 7 edge-disjoint triangles, each representing a line in the poem.}

\end{table}
\end{example}

\bigskip

\begin{example}
Relaxed Steiner triple poem: 
\\  

\textbf{A Pause in the Rain}
\\
\\
\indent \textbf{Rain} fell; the \textbf{man} sipped \textbf{coffee}.

\textbf{Cold} \textbf{rain}; all \textbf{wet}!

In \textbf{cold} the \textbf{man} watched a \textbf{bird}.

The \textbf{wet}  \textbf{man} paused in \textbf{wonder}.

With \textbf{coffee}, in \textbf{cold} \textbf{wonder}.

The \textbf{bird} chirped in the \textbf{rain}, in \textbf{wonder}.

The \textbf{wet} \textbf{bird} alighted as \textbf{coffee} steamed!

\vspace{0.4 cm}

Keywords: \{Rain, Man, Coffee, Cold, Wet, Bird, Wonder \}

\vspace{0.4 cm}

\textbf{Discussion of the poem:}  This poem, written by the first author (De Vas Gunasekara), emerged from an interest in exploring the combinatorial structure of a Steiner triple system through a poetic lens. While the initial motivation was mathematical - she was using the structure of $\mathrm{STS}(7)$ as a formula, the creative process evolved into a personal reflection. Drawing on her love for nature and photography, De Vas Gunasekara imagined a quiet, rainy scene that resonates with her own enjoyment of peaceful moments, like sipping coffee on a wet evening.
To choose the keywords, she focused on elements that are naturally interconnected and rich in imagery: rain, cold, wet, man, coffee, bird, and wonder. The use of man instead of a direct self-reference allowed her to create distance, framing the moment as something observed rather than experienced firsthand.
The poem captures a fleeting scene of quiet reflection, where a man, caught in the rain with a cup of coffee, observes a bird and finds himself drawn into a moment of shared stillness. The combinatorial structure supports this calmness, while the simple language deepens its emotional resonance.

\vspace{0.4 cm}

Note that in both of these examples, in addition to following a Steiner Triple System, there is an additional property: any two lines in the poem share exactly one common keyword. This means the structure also corresponds to that of a so-called \emph{Fano plane}. This is, in fact, due to all \(\mathrm{STS}(7)\) being Fano planes.

\end{example}

\subsection{\texorpdfstring{$\mathrm{STS}(9)$}{STS(9)}  - Poems with nine Keywords}

$\mathrm{STS}(9)$ is a Steiner Triple System of order 9. A poem that conforms to this structure contains nine distinct keywords. Since each pair of keywords must appear together exactly once, the poem consists of $\frac{9(9 - 1)}{6} = 12$ lines, with exactly three keywords per line. Consequently, each keyword appears in precisely four different lines.

Below, we present our constructions of pure Steiner Triple poem, a relaxed version, and a resolvable version. We also illustrate the underlying graph representations of these poems.

\begin{example}
 Relaxed Steiner triple poem: 
 \\
 \\
\indent \textbf{Footprints on a Snowy Evening}

\indent (after  ``Stopping by the Woods on a Snowy Evening" by   Robert Frost \cite{Frost1923}) 
\\
\\
\indent \textbf{Dark} \textbf{woods}, faded \textbf{footprints}.

\textbf{Woods}, \textbf{lake} and a \textbf{man}.

\textbf{Man} rambles, \textbf{dark} \textbf{wind} whispers.

The \textbf{woods} are cold under \textbf{snow} and \textbf{wind}.

\textbf{Footprints} silent, lake \textbf{frozen}, \textbf{snow} falls.

\textbf{Evening} is \textbf{dark} by the \textbf{lake}.

\textbf{Evening} deepens, \textbf{snow} brushes the \textbf{man}.

A \textbf{horse} drifts through \textbf{evening} \textbf{woods}.

The \textbf{horse} shivers, \textbf{lake} still, \textbf{wind} howls.

\textbf{Evening} falls; lingering  \textbf{footprints}, \textbf{wind} murmurs gently.

\textbf{Dark} \textbf{horse} blows, \textbf{snow} drifts.

\textbf{Man} walks, \textbf{horse} follows-\textbf{footprints} remain.

\vspace{0.4 cm}

\textbf{Discussion of the poem:} ``Footprints on a Snowy Evening" is a poem written by the first author (De Vas Gunasekara) as part of her ongoing exploration of combinatorial structures in poetry. This piece began with a formal constraint using nine interconnected words arranged according to a combinatorial pattern of $\mathrm{STS}(9)$. Initially uncertain about the content, she turned to inspiration from existing literature, particularly Robert Frost’s ``Stopping by Woods on a Snowy Evening", one of her longtime favourites.
The poem acts as both a structural experiment and a subtle homage. After selecting nine thematically linked words to evoke a quiet, wintry landscape, De Vas Gunasekara constructed the scene using repetition of the keywords to deepen atmosphere and emotional resonance.
The mathematical structure ensures that each word pairs with every other exactly once, adding a hidden layer of cohesion.
Following feedback from the second author, the poem was redrafted to strengthen narrative clarity and poetic flow. The deliberate repetition within the design creates an intensifying effect, transforming a static winter moment into something softly evocative.

To visualize the structure graph-theoretically, we represent each keyword as a vertex and construct the complete graph on nine vertices, denoted $K_9$. Each line of the poem corresponds to a triangle in this graph. Therefore, the poem forms a triangle decomposition of $K_9$, in which the twelve triangles collectively cover all the edges. Each triangle, representing a line of the poem, contains exactly three keywords.

Keywords (vertex set): \\ \indent \{Dark, Woods, Man, Lake, Footprints, Wind, Snow, Evening, Horse  \}

\begin{figure}[h!]
\centering
\begin{tikzpicture}[scale=4, every node/.style={circle, draw, minimum size=1.95cm, font=\small}]
  % Define the 9 nodes in a circular layout
  \node (Dark)       at (90:1)    {Dark};
  \node (Woods)      at (50:1)    {Woods};
  \node (Man)    at (10:1)    {Man};
  \node (Lake)       at (-30:1)   {Lake};
  \node (Footprints)  at (-70:1)   {Footprints};
  \node (Wind)       at (-110:1)  {Wind};
  \node (Snow)       at (-150:1)  {Snow};
  \node (Evening)    at (-190:1)  {Evening};
  \node (Horse)      at (130:1)   {Horse};

  % 12 edge-disjoint triangles from your poem
  \draw[thick, red]             (Dark) -- (Woods) -- (Man) -- (Dark) --cycle;     % Line 1
  \draw[thick, blue]            (Woods) -- (Lake) -- (Footprints) -- (Woods) -- cycle;     % Line 2
  \draw[thick, green!70!black]  (Footprints) -- (Dark) -- (Wind) -- (Footprints)--cycle;      % Line 3
  \draw[thick, orange]          (Man) -- (Lake) -- (Snow) -- (Man) -- cycle;        % Line 4
  \draw[thick, magenta]         (Woods) -- (Snow) -- (Wind) -- (Woods)--cycle;          % Line 5
  \draw[thick, teal]            (Evening) -- (Dark) -- (Lake) --(Evening)-- cycle;        % Line 6
  \draw[thick, violet]          (Horse) -- (Evening) -- (Woods) --(Horse)-- cycle;      % Line 7
  \draw[thick, cyan]            (Horse) -- (Lake) -- (Wind) -- (Horse)--cycle;          % Line 8
  \draw[thick, brown]           (Evening) -- (Man) -- (Wind) --(Evening)-- cycle;     % Line 9
  \draw[thick, lime!80!black]   (Evening) -- (Snow) -- (Footprints) --(Evening)-- cycle;   % Line 10
  \draw[thick, gray]            (Man) -- (Footprints) -- (Horse) --(Man)-- cycle;  % Line 11
  \draw[thick, pink]            (Dark) -- (Horse) -- (Snow) --(Dark)-- cycle;          % Line 12

\end{tikzpicture}
% \end{center}
\caption{Triangle decomposition of the complete graph $K_9$ on nine keywords from the poem ``Footprints on a Snowy Evening". Each triangle corresponds to one line in the relaxed Steiner triple poem.}
\label{fig:sts7-poem2}
\end{figure}
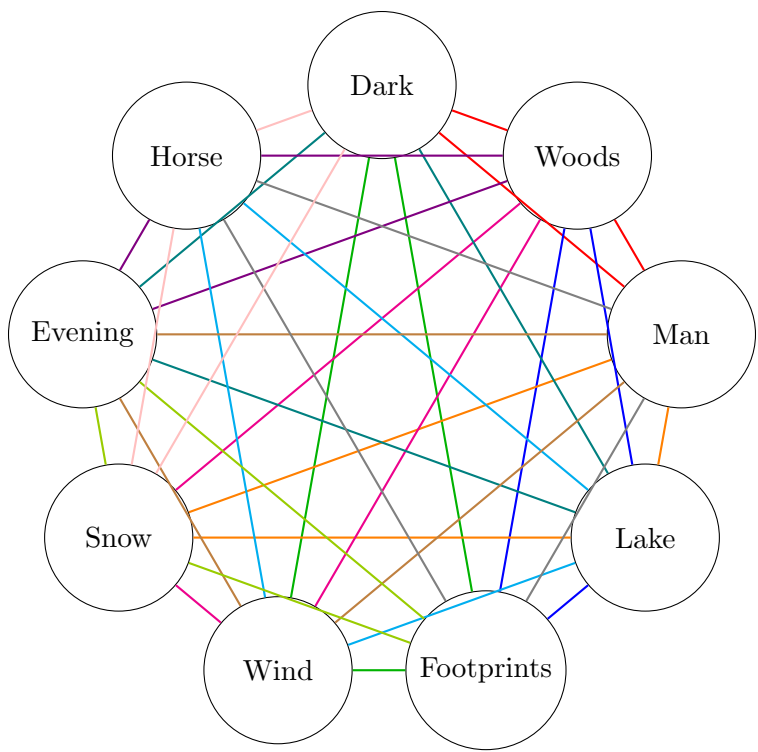

\begin{table}[H]
\centering
\begin{tabular}{@{}lll@{}}
\toprule
\textbf{Line \#} & \textbf{Triangle (Color)} & \textbf{Vertices} \\
\midrule
1  & Red        & \{Dark, Woods, Man\} \\
2  & Blue       & \{Woods, Lake, Footprints\} \\
3  & Green      & \{Footprints, Dark, Wind\} \\
4  & Orange     & \{Man, Lake, Snow\} \\
5  & Magenta    & \{Woods, Snow, Wind\} \\
6  & Teal       & \{Evening, Dark, Lake\} \\
7  & Violet     & \{Horse, Evening, Woods\} \\
8  & Cyan       & \{Horse, Lake, Wind\} \\
9  & Brown      & \{Evening, Man, Wind\} \\
10 & Lime       & \{Evening, Snow, Footprints\} \\
11 & Gray       & \{Man, Footprints, Horse\} \\
12 & Pink       & \{Dark, Horse, Snow\} \\
\bottomrule
\end{tabular}
\caption{The 36 edges of \( K_9 \) are partitioned into 12 edge-disjoint triangles, each representing a line in the poem.}
\label{tab:triangle-k9}
\end{table}

\end{example}

%\vspace{0.4 cm}

\begin{example}
\label{E:STS9 pure}
 Pure Steiner triple poem: 
 \\
 
 \textbf{Wordstorm }
 
 (after ``Night Feed” by Debbie Lim \cite{Lim2021})
\\

Face! Look, poem! 

Poem! Fix lines! 

Lines! Look helpless! 

Helpless! Poem river! 

River! Warm lines!

Lines! Call face! 

Face! Warm helpless! 

Helpless! Fix call! 

Call! Look river!

River! Fix face! 

Fix! Look warm! 

Warm! Call poem.
\end{example}

\vspace{0.3 cm}
Keywords: \{Face, Look, Poem, Fix, Lines, Helpless, River, Warm, Call\}

\vspace{0.3 cm}
\textbf{Discussion of the poem:} It occurred to the second author (Lo) as she chose keywords from Lim's poem, that words which could double as both verbs and nouns (e.g. lines, call, fix) were particularly useful. Lo experienced an unexpected pleasure in the act of determining the poetic constraints (by drawing the balls-and-triangles graph and noting down the word groupings) - it felt very calming. The process drew her into a meditative space, comparable to that experienced when doing Zentangle art. However, when it came to composing the poem, Lo experienced frustration. She wanted to write with expressive intent, preferably about the experience of breastfeeding her own children. This was much more challenging within the given constraints of a pure Steiner triple poem (as she was restricted to using keywords from the poem, without any filler words). She had to let go of her original intent and allow the constraints to have greater power over the construction of the poem. She quickly switched to using imperatives because this was one way to create a semblance of meaning. She decided to add an extra constraint to create a sense of continuity: the last word of each line repeats as the first word of the next line. The resulting poem is much more Oulipo (in its resistance to expressivism) than Lo originally intended.

\begin{example}
 Resolvable Steiner triple poem: 
 \\
 
 \textbf{Things We Cannot Keep }
 
 (after ``Night Feed”   by Debbie Lim \cite{Lim2021}) 
\\

\textbf{Helpless} words, for this \textbf{warm} baby \textbf{face} 

at my breast: a \textbf{poem} I cannot \textbf{fix} into \textbf{lines}.

I \textbf{call} to the \textbf{river}: ``Stop! \textbf{Look!}''

\vspace{1em}

\textbf{Helpless}, the \textbf{poem} makes its \textbf{call}: 

\textbf{warm} words, to \textbf{fix} the coldness of the \textbf{river}; 

a soft \textbf{face}, urgent \textbf{lines}, an unguarded \textbf{look}.

\vspace{1em}

\textbf{Helpless}, a mother in a \textbf{fix}, where will I \textbf{look} 

if \textbf{warm} \textbf{lines} cannot \textbf{call} 

that baby \textbf{face} back? I throw the \textbf{poem} in the \textbf{river}.

\vspace{1em}

\textbf{Helpless} \textbf{lines} surface like \textbf{river} cormorants 

with \textbf{warm} fish for an impossible \textbf{poem}; don’t \textbf{look} back; 

\textbf{face} what you cannot \textbf{fix}; let the \textbf{call}, fall. 
\end{example}

\vspace{0.3 cm}
Keywords: \{Face, Look, Poem, Fix, Lines, Helpless, River, Call, Warm\}

\vspace{0.3 cm}
\textbf{Discussion of the poem:} This poem, by the second author (Lo), uses the same keywords as ``Wordstorm" (Example~\ref{E:STS9 pure}), but with two differences in constraints: first, the resolvable constraint, which dictated that all keywords be contained in each three-line stanza; second, that filler words were allowed. The ability to use filler words gave Lo much more scope to achieve her expressivist intent for this poem, which was to write about her own experience breastfeeding her children. Lo also wanted to amplify the theme of Lim's original poem, ``Night Feed", which is about the passing of time and the acceptance of loss. Lo is unsure of how well the resultant poem stands on its own; however, the constraints of patterned repetition did allow for a poetic exploration and intensification of existing themes from Lim's original poem.

\section{Discussion}

The first author comes from a mathematical background, while the second author comes from a literary background. The experience of constructing poetry under mathematical constraints affected the two authors in contrasting ways. The first author, De Vas Gunasekara often began with a set of interconnected words or drew inspiration from an existing poem, then applied a mathematical structure to generate flow and meaning. For her, the process felt like following a recipe - systematic and guided. The constraints provided by the mathematical design helped her shape the poem, offering a framework that supported the creative act.

By contrast, the second author, Lo (who has a more literary background), described a more complex and emotionally layered response. Her writing process began with free writing to determine a subject choice: the free writing helped her to figure out what topic, theme or issue she was drawn to and why. She would \textit{then} consider how her desire to write about a chosen topic could interact with the given form. During this stage of writing, she would experience the constraints limiting her expressive freedom. Often the sense of limitation would be accompanied by a paradoxical sense of liberation. Lo reflected:

\begin{quote}
I feel like I am making choices that will limit what I can do in the poem … but that feels good, because writing is all about making choices. Now I have fewer choices, which gives me a strange sense of freedom. Freedom within limits, maybe that is part of what poetry is?    
\end{quote}

In relation to a different poem, Lo noted:

\begin{quote}
How secure and calm the process makes me feel … even though I am doing this in a timed setting; I am drawn into the process, I can suspend time, even with stress pressing in – all the rules help (weirdly). Do this, then do that. Zen-tangle poetry? (I think of Adolin playing towers with Yanagawn as Odium’s armies mount their invasion) ( \cite{Sanderson2024} Sanderson, 2024).
I am doing this a line at a time, writing a line and letting meaning, flow and repetition determine which set I pick next.    
\end{quote}

At other times, the second author experienced frustration with the imposed structure, and felt dissatisfied with the poetic outcome. 
However, even this dissatisfaction could be framed in a productive way. From a poetry-practice perspective, working within the tight constraints of Steiner Triple Systems could be a useful exercise for poets exploring the tension between formal constraints and expressivist intent. At what point do the constraints become too restrictive? When and where could a poet break the rules? These are really interesting questions and could prove useful in generating self-awareness within a poetry workshop setting. 

When the writing experiences of the first and second authors  are compared, it is clear that mathematical frameworks can both guide and challenge creative processes, shaping each author's engagement with poetry in uniquely personal ways.

\section{Conclusion}

This work opens several avenues for future exploration. In this paper, we focused on poetic constructions based on Steiner triple systems of orders seven and nine. However, Steiner triple systems exist for all positive integers \( u \equiv 1,3 \pmod{6} \), allowing for the creation of similar poetic structures using a broader range of word counts. More generally, one can extend this framework to block designs with block sizes greater than three, enabling lines or stanzas that incorporate four or more key words. Other combinatorial structures, such as Latin squares or finite geometries, may also be adapted to guide poetic composition, offering rich terrain for the development of new poetic forms grounded in mathematical logic.

Incorporating combinatorial design into poetic activities offers a novel and interdisciplinary approach to mathematics education. For students in the humanities, this method provides an accessible entry point into mathematical thinking, allowing them to engage with concepts such as structure, symmetry, and constraint in a creative and meaningful way. By working within mathematically defined patterns such as those imposed by Steiner triple systems, students develop intuition for combinatorial reasoning while actively constructing poetic texts. In doing so, it can demystify abstract mathematical ideas and promote confidence in working with them.

Conversely, for students with a background in mathematics, this approach offers a valuable opportunity to experience mathematical structures not only as formal objects but also as generative frameworks for artistic expression. The interplay between free writing and rule-based composition highlights the creative potential of constraint and offers insight, not only into the aesthetics of structure, but also into tensions between content and form. In this way, the project fosters a bidirectional appreciation: students from literary disciplines gain mathematical insight, while students from mathematical disciplines cultivate creative flexibility and communication through interdisciplinary exploration.


\begin{thebibliography}{9}

\bibitem{PoetryFoundationForm}
H. Amos, M. Queeney, and R. E. Shoemaker,
\textit{Poetry and Form}, Poetry Foundation,
\url{https://www.poetryfoundation.org/collections/159875/poetry-and-form} 



\bibitem{Barriere2017}
    L. Barrière,
    Combinatorics in the Art of the Twentieth Century,
    {\it Proceedings of Bridges 2017: Mathematics, Art, Music, Architecture, Education, Culture},
    (2017), 321--328. \\
    \url{https://archive.bridgesmathart.org/2017/bridges2017-321.pdf}

\bibitem{Baumgaertner1990}
J. P. Baumgartner, 
 Poetry, 
{Harcourt Brace Jovanovich, New York}, 1990.

    


 \bibitem{ColRos}
    C.J. Colbourn and A. Rosa,
    Triple Systems,
    {\it Clarendon Press, Oxford} (1999). 

    
    

    \bibitem{Frost1923}
R. Frost, Stopping by Woods on a Snowy Evening, 
originally published in \textit{New Hampshire} (1923). 
Available at: \url{https://www.poetryfoundation.org/poems/42891/stopping-by-woods-on-a-snowy-evening}

    

    \bibitem{Glaz2011}
S. Glaz, 
Poetry inspired by mathematics: A brief journey through history, 
{\it J. Mathematics Arts} {\bf 5}(4) (2011), 171--183.

\bibitem{Hart2023}
    S. Hart, 
    Once Upon a Prime: The Wondrous Connections Between Mathematics and Literature, 
    {\it Flatiron Books, New York}, 2023.

   \bibitem{Hass2017}
   R. Hass, 
   A Little Book on Form: An Exploration into the Formal Imagination of Poetry, 
   {\it Ecco (HarperCollins), New York}, 2017.
   
   

\bibitem{HetheringtonAtherton2020}
P. Hetherington and C. Atherton,
{\it Prose Poetry: An Introduction}, Princeton U.P., Princton, 2020.

    \bibitem{KaraaliLesser2020}
G. Karaali and L. M. Lesser, 
Mathematics and Poetry: Arts of the Heart,
{\it Handbook of the Mathematics of the Arts and Sciences}, Springer, Cham (2020), 1--13.

\bibitem{kirkman1847}
    T.P. Kirkman,
    On a problem in combinations,
    {\it Cambridge and Dublin Math. J.} {\bf 2} (1847), 191--204.

    

\bibitem{Lim2021}
D. Lim, Night Feed, In E. Kurz, S. King and C. Delahunty (Eds.), {\it What We Carry} (p. 156), Recent Work Press, Canberra, 2021. 

\bibitem{Lu1965}
J. X. Lu, Collected Works of Lu Jiaxi on Combinatorial Designs,
{\it Inner Mongolia People’s Press}, (1965). 


    \bibitem{May2021}
    D. May,
    Poems Structured by Mathematics,
    {\it Handbook of the Mathematics of the Arts and Sciences},
    (2021), 1045--1092. \\
    \url{https://link.springer.com/referenceworkentry/10.1007/978-3-319-57072-3_113}


    \bibitem{Mirsky1979}
    L. Mirsky,
    Book Review: Combinatorics with emphasis on the theory of graphs,
    {\it Bull. Amer. Math. Soc. (N.S.)} {\bf 2(1)} (1979), 380--388.


    \bibitem{RayChaudhuriWilson1971}
    D. K. Ray-Chaudhuri and R. M. Wilson,
    Solution of Kirkman’s schoolgirl problem ,
    {\it Combinatorics, Proc. Sympos. Pure Math. Vol. XIX}  (1971), 187--203.

    \bibitem {Sanderson2024}
    B. Sanderson,
    \textit{Wind and Truth}, Gollancz, 2024. 

    
\end{thebibliography}
\end{document}